\documentstyle{article}
\begin{document}
\begin{center}
{\bf{\LARGE { Edge Theorem for Multivariable Systems \footnote[1]
{ Supported by National Natural Science Foundation of
China(69925307). Email: longwang@mech.pku.edu.cn} }}}
\end{center}

\vskip 0.6cm \centerline{ Long Wang$^1$ \hspace{2cm} Zhizhen
Wang$^1$ \hspace{2cm} Lin Zhang$^2$ \hspace{2cm}Wensheng Yu$^3$}
\vskip 6pt \centerline{{$^1$ Center for Systems and Control,
Department of Mechanics, Peking University, Beijing 100871,
China}} \centerline{{$^2$ Computer Science Department, Naval
Postgraduate School, Monterey, CA93943, USA }} \centerline{{$^3$
Institute of Automation, Chinese Academy of Sciences, Beijing
100080, China }}

\vskip 0.6cm

\vskip 6pt
\begin{minipage}[t]{14cm}{
{\bf Abstract: } This paper studies robustness of multivariable
systems with parametric uncertainties, and establishes a
multivariable version of Edge Theorem. An illustrative example is
presented.

\vskip 4pt

 }
\end{minipage}

\section{\bf{Introduction}}

\hspace{1cm} Motivated by the seminal theorem of Kharitonov on robust stability of interval
polynomials\cite{Khar, Khar1}, a number of papers on robustness analysis of uncertain
systems have been published in the past few years\cite{Holl, Bart, Fu, wang1,
wang2, Barmish, Chap, wang3}. Kharitonov's theorem states that the Hurwitz stability of
the real (or complex) interval polynomial family can be guaranteed by the Hurwitz
stability of four (or eight) prescribed critical vertex polynomials in this family.
This result is significant since it reduces checking stability of infinitely many
polynomials to checking stability of finitely many polynomials, and the number of
critical vertex polynomials need to be checked is independent of the order of the
polynomial family. An important extension of Kharitonov's theorem is the edge theorem
discovered by Bartlett, Hollot and Huang\cite{Bart}. The edge theorem states that
the stability of a polytope of polynomials can be guaranteed by the stability of its
one-dimensional exposed edge polynomials. The significance of the edge theorem is
that it allows some (affine) dependency among polynomial coefficients, and applies to
more general stability regions, e.g., unit circle, left sector, shifted half plane,
hyperbola region, etc. When the dependency among polynomial coefficients is nonlinear,
however, Ackermann shows that checking a subset of a polynomial family generally can
not guarantee the stability of the entire family\cite{Ack1, Ack2}.

For Hurwitz stability of interval matrices, Bialas 'proved' that in order to guarantee robust stability,
it suffices to check all vertex matrices\cite{Bialas}. Later, it was shown by Barmish that Bialas' result
was incorrect\cite{Barmish1}. Kokame and Mori eastblished a Kharitonov-like result on robust Hurwitz
stability of interval polynomial matrices\cite{Mori}, and Kamal and Dahleh established some robust stability
criteria for MIMO systems with fixed controllers and uncertain plants\cite{kama}.

In this paper, we will study robustness of a class of MIMO systems
with their transfer function matrices described by

\begin{equation}\label{1}
\begin{array}{l}
{\cal F}(s)=\left\{\left(\begin{array}{lll}
a_{11}(s)&\dots&a_{1n}(s)\\
\dots&\dots&\dots\\
a_{n1}(s)&\dots&a_{nn}(s)
\end{array}\right): \  \ a_{ij}(s)\in {\cal A}_{ij}(s)\right\}\\
{\cal A}_{ij}(s){=}\mbox{conv}\left\{b1_{ij}(s),
\dots,bm_{ij}(s)\right\}
\end{array}
\end{equation}
where $m$ is a given positive integer.

\section{\bf{Preliminaries}}
{

{\bf Definition 1}\quad  A polynomial matrix is a matrix with all its entries being polynomials.

{\bf Definition 2}\quad  Suppose $D$ is a simply-connected region in the complex plane. If all the roots of the determinant of
a polynomial matrix lie within $D$, then this polynomial matrix is
called $D$-stable. A set of polynomial matrices is called robustly $D$-stable, if
every polynomial matrix in this set is $D$-stable.

{\bf Definition 3}\quad Suppose $f_1(s),\dots,f_m(s)$ are $m$ given polynomials, the set
$$
\left\{\sum_{i=1}^m \lambda_i f_i(s):\ \ \lambda_i\geq 0, \ \ \sum_{i=1}^m \lambda_i=1 \right\}
$$
is called the polynomial polytope generated by $f_1(s),\dots,f_m(s)$, denoted as $\mbox{conv}\left\{f_1(s),
\dots,f_m(s)\right\}$.

{\bf Definition 4}\quad The polynomial $f(s)=a_0+a_1 s+\dots+a_n s^n$, with $a_i\in[a_i^L,a_i^U]$
is called an interval polynomial.

{\bf Definition 5}\quad $S_n$ is the set of all bijections from $\{1,\dots,n\}$ to
$\{1,\dots,n\}$.

{\bf Definition 6}\quad The vertex set and edge set of ${\cal A}_{ij}(s)$ are
$$
K_{ij}(s)=\{b1_{ij}(s) , \dots, bm_{ij}(s) \}
$$
$$
E_{ij}(s)=\{ \lambda br_{ij}(s)+(1-\lambda)bt_{ij}(s),  \quad \lambda \in [0, 1],  \quad r,t\in\{1, \dots, m \}\}
$$
respectively.

{\bf Definition 7}\quad
\begin{equation}
{\cal F}_E(s)=\bigcup_{\sigma\in S_n }\left \{(a_{ij}(s))_{n \times n}:
a_{ij}(s)\left\{
\begin{array}{l}
\in E_{ij}(s)  \mbox{  if   } i=\sigma(j) \\
 \in K_{ij}(s) \mbox{  if   } i\not=\sigma(j)
\end{array}\right.
\right \}
\end{equation}

{\bf Lemma 1}\quad (Edge Theorem\cite{Bart}) Suppose $\Gamma \subset
{\cal C}$ is a simply-connected region, $\Omega$ is a polynomial polytope without degree dropping.
Then, $\Omega$ is $\Gamma$-stable if and only if all the edges of $\Omega$ are $\Gamma$-stable.

{\bf Lemma 2}\quad Suppose $A(s)$ is a given ${n\times (n-1)}$ polynomial
matrix. Then
$$
\begin{array}{l}
\left\{\left( \begin{array}{c}
a_{11}(s)\\
\vdots\\
a_{n1}(s)
\end{array}
A(s)\right): \ \
\begin{array}{c}
a_{i1}(s)\in{\cal A}_{i1}(s),\\
 i=1,\dots,n
 \end{array}
\right\}
 \mbox{  is robustly $D$-stable  }\\
\Leftrightarrow \mbox{ for all $i=1,\dots,n$,}
\left\{\left( \begin{array}{c}
a_{11}(s)\\
\vdots\\
a_{n1}(s)
\end{array}
A(s)\right):\ \ \begin{array}{ll}
a_{l1}(s)\in K_{lj}(s)& l\not=i\\
a_{l1}(s)\in E_{lj}(s)& l=i
\end{array}\right\} \mbox{  is robustly $D$-stable.  }
\end{array}
$$

Proof:\ \ Necessity is obvious, since the later is a subset of the former.

Sufficiency: For any $a_{i1}(s)\in {\cal A}_{i1}(s)$, the corresponding matrix is
$$
T(s)=\left( \begin{array}{c}
a_{11}(s)\\
\vdots\\
a_{n1}(s)
\end{array}
A(s)\right)
$$
By Laplace Formula, we can expand the determinant of $T(s)$ along its first column. Then,
by convexity and by Lemma 1, we know that
$T(s)\mbox{ is robustly $D$-stable}$.

{\bf Lemma 3}\quad Suppose $B(s)$ is a given ${(n-1)\times n}$ polynomial matrix.
$\ast$ stands for fixed entries in a matrix. Then
$$
\begin{array}{l}
\left\{\left( \begin{array}{c}
\begin{array}{lllll}
\ast & a_{1i}(s) & \ast & a_{1j}(s)& \ast \\
\end{array}\\
B(s)
\end{array}
\right):
\begin{array}{l}
a_{1i}(s) \in  {\cal A}_{1i}(s)\\
a_{1j}(s) \in  {\cal A}_{1j}(s)
\end{array}\right\}
 \mbox{  is robustly $D$-stable }\\

\Leftrightarrow
\left\{\left(\begin{array}{c}
 \begin{array}{lllll}
\ast & a_{1i}(s) & \ast & a_{1j}(s)& \ast
\end{array}\\
B(s)
\end{array}\right):
\begin{array}{c}
a_{1i}(s)\times a_{1j}(s)\in \\
 \left(K_{1i}(s)\times E_{1j}(s)\right)\cup \left(E_{1i}(s)\times K_{1j}(s)\right)
 \end{array}\right\}
 \mbox{  is robustly $D$-stable.  }
\end{array}
$$

Proof:\ \ the proof is analogous to the proof of Lemma 2, except that the Laplace expansion is
carried out along the row instead of the column.

\section{ Main Results }

{\bf Theorem 1}\quad ${\cal F}(s)$ is robustly $D$-stable if and only if ${\cal F}_E(s)$
is robustly $D$-stable.

Proof:\ \ Necessity is obvious. To prove sufficiency, we first note that interchanging any two rows (or columns) does not affect
the stability of a polynomial matrix (it only changes the sign of the determinant). By Lemma 2
$$
\begin{array}{c}
{\cal F}(s)\mbox{ is robustly $D$ stable }\\
\Leftrightarrow \mbox{ for all $i=1,\dots,n$,}
\left\{\left( \begin{array}{c}
a_{11}(s)\\
\vdots\\
a_{n1}(s)
\end{array}
A(s)\right):\ \ \begin{array}{ll}
a_{l1}(s)\in K_{lj}(s)& l\not=i\\
a_{l1}(s)\in E_{lj}(s)& l=i
\end{array}\right\} \mbox{  is robustly $D$ stable.  }\\
\Leftrightarrow \mbox{ for all }
\left\{\begin{array}{l}
i_1=1,\dots,n \\
i_2=1,\dots,n
\end{array}\right.
\left\{\left( \begin{array}{cc}
a_{11}(s)&a_{12}(s)\\
\vdots&\vdots\\
a_{n1}(s)&a_{n2}(s)
\end{array}
A_1(s)\right):\ \
\begin{array}{l}
\left\{\begin{array}{ll}
a_{l1}(s)\in K_{lj}(s)& l\not=i_1\\
a_{l1}(s)\in E_{lj}(s)& l=i_1
\end{array}\right.\\
\left\{\begin{array}{ll}
a_{l2}(s)\in K_{lj}(s)& l\not=i_2\\
a_{l2}(s)\in E_{lj}(s)& l=i_2
\end{array}\right.
\end{array}
\right\} \\
\mbox{  is robustly $D$ stable.  }
\end{array}
 $$
where $A_1(s)$ is the corresponding ${n\times (n-2)}$ polynomial matrix. This last equivalence is based on Lemma 2 and the fact
 that interchanging two columns does not change the stability  of a polynomial matrix.
Repeating the process above, let $Y_n$ denote the set of all mappings from $\{1,\dots,n\}$ to
 $\{1,\dots,n\}$, then
$$
\begin{array}{l}
{\cal F}(s)\mbox{ is robustly $D$ stable}\\
\Leftrightarrow \mbox{ for all $\eta\in Y_n$,}
\left\{\left( a_{ij}(s)\right):\ \ \begin{array}{ll}
a_{ij}(s)\in K_{ij}(s)& i\not=\eta(j)\\
a_{ij}(s)\in E_{ij}(s)& i=\eta(j)
\end{array}\right\} \mbox{  is robustly $D$ stable.  }
\end{array}
$$
If there exists an $\eta\in Y_n$ such that
$\eta(i_1)=\eta(i_2)=k$, then the corresponding matrix $F(s)=(a_{ij}(s))$ satisfies
$$
\begin{array}{l}
a_{i_1k}(s)\in E_{i_1k}(s)\\
a_{i_2k}(s)\in E_{i_2k}(s)
\end{array}
$$
Applying Lemma 2 to column $k$ of
$F(s)$, we have
$$
\begin{array}{l}
{\cal F}(s)\mbox{ is robustly $D$ stable}\\
\Leftrightarrow \mbox{ for all $\sigma\in S_n$,}
\left\{\left( a_{ij}(s)\right):\ \ \begin{array}{ll}
a_{ij}(s)\in K_{ij}(s)& i\not=\sigma(j)\\
a_{ij}(s)\in E_{ij}(s)& i=\sigma(j)
\end{array}\right\} \mbox{  is robustly $D$ stable.  }\\
\Leftrightarrow {\cal F}_E(s) \mbox{  is robustly $D$ stable.  }
\end{array}
$$

\section{Interval Model}

Interval model, as a simple and effective approximation of uncertain
systems, has been the subject of study in robustness analysis for
a long time. In a similar vein, we consider the Hurwitz stability
of the following uncertain system.
\begin{equation}
\begin{array}{l}
  {\cal G}(s)=\left\{(c_{ij}(s)):c_{ij}(s)\in{\cal C}_{ij}(s)\right\}\\
  {\cal C}_{ij}(s)\mbox{ are interval polynomials}
  \end{array}
\end{equation}

{\bf Definition 8}\quad For the interval polynomial ${\cal C}_{ij}(s)=
\{\sum_{l=0}^n q_l(ij) s^l, \quad q_l(ij) \in [\underline q_l(ij), \overline q_l(ij)]
\}$, its Kharitonov vertex set and Kharitonov edge set are defined respectively as
$$
\begin{array}{c}
K_{ij}^I(s)=\{c_k^1(s), c_k^2(s), c_k^3(s), c_k^4(s)\}\\
E_{ij}^I(s)=\{\lambda c_k^r(s)+(1-\lambda) c_k^t(s), \quad
(r,t)\in \{(1,2), (2,4), (4,3), (3,1)\}, \quad  \lambda \in [0,1]\}
\end{array}
$$
where
$$
  \begin{array}{ll}
  c_k^1(s)=\underline q_0(ij)+\underline q_1(ij) s+\overline q_2(ij) s^2+
  \overline q_3(ij) s^3+\dots \quad & \quad
  c_k^2(s)=\underline q_0(ij)+\overline q_1(ij) s+\overline q_2(ij) s^2+
  \underline q_3(ij) s^3+\dots \\
  c_k^3(s)=\overline q_0(ij)+\underline q_1(ij) s+\underline q_2(ij) s^2+
  \overline q_3 s^3+\dots \quad & \quad
  c_k^4(s)=\overline q_0(ij)+\overline q_1(ij) s+\underline q_2(ij)v s^2+
  \underline q_3(ij) s^3+\dots
  \end{array}
  $$

{\bf Definition 9}\quad
\begin{equation}
{\cal G}_E(s)=\bigcup_{\sigma\in S_n }\left \{(c_{ij}(s))_{n \times n}:
c_{ij}(s)\left\{
\begin{array}{l}
\in E_{ij}^I(s)  \mbox{  if   } i=\sigma(j) \\
 \in K_{ij}^I(s) \mbox{  if   } i\not=\sigma(j)
\end{array}\right.
\right \}
\end{equation}

{\bf Lemma 4}\quad (Box Theorem\cite{kama}) Suppose
$\Delta(s)=\{\delta(s,p)=F_1(s) P_1(s)+\dots+F_m(s) P_m(s)\}$,
  $P_i(s)$ is an interval polynomial, $F_i(s)$ is a given fixed polynomial, $i=1, \dots, m$.
 And suppose $\Delta(s)$ is degree-invariant. Then, $\Delta(s)$ is Hurwitz stable if and only if
 $\Delta_E (s)$ is Hurwitz stable, where
$\Delta_E (s)=\cup_{l=1}^m \{\sum_{i=1}^{l-1}F_i (s)K_{P_i}^0 (s)+
F_l(s) E_{P_l}^0(s)+\sum_{i=l+1}^m F_i (s)K_{P_i}^0(s)\}$
(let $\sum_{i=r}^t f_i=0$, if $r>t$).

By resort to the Box Theorem, and following a similar line of arguments as in the proof of Theorem 1,
we can get some analogous stability verification results for the interval model.

{\bf Lemma 5}\quad Suppose $A(s)$ is a given ${n\times (n-1)}$ polynomial matrix.
Then
$$
\begin{array}{l}
\left\{\left( \begin{array}{c}
c_{11}(s)\\
\vdots\\
c_{n1}(s)
\end{array}
A(s)\right): \ \
\begin{array}{c}
c_{i1}(s)\in{\cal G}_{i1}(s),\\
 i=1,\dots,n
 \end{array}
\right\}
 \mbox{  is robustly Hurwitz stable.  }\\
\Leftrightarrow \mbox{ for all $i=1,\dots,n$,}
\left\{\left( \begin{array}{c}
c_{11}(s)\\
\vdots\\
c_{n1}(s)
\end{array}
A(s)\right):\ \ \begin{array}{ll}
c_{l1}(s)\in K_{lj}^I(s)& l\not=i\\
c_{l1}(s)\in E_{lj}^I(s)& l=i
\end{array}\right\} \mbox{  is robustly Hurwitz stable.  }
\end{array}
$$

{\bf Lemma 6}\quad Suppose $B(s)$ is a given ${(n-1)\times n}$ polynomial matrix.
$\ast$ stands for fixed entries in a matrix. Then
$$
\begin{array}{l}
\left\{\left( \begin{array}{c}
\begin{array}{lllll}
\ast & c_{1i}(s) & \ast & c_{1j}(s)& \ast \\
\end{array}\\
B(s)
\end{array}
\right):
\begin{array}{l}
c_{1i}(s) \in  {\cal C}_{1i}(s)\\
c_{1j}(s) \in  {\cal C}_{1j}(s)
\end{array}\right\}
 \mbox{  is robustly Hurwitz stable }\\

\Leftrightarrow
\left\{\left(\begin{array}{c}
 \begin{array}{lllll}
\ast & c_{1i}(s) & \ast & c_{1j}(s)& \ast
\end{array}\\
B(s)
\end{array}\right):
\begin{array}{c}
c_{1i}(s)\times c_{1j}(s)\in \\
 \left(K_{1i}^I(s)\times E_{1j}^I(s)\right)\cup \left(E_{1i}^I(s)\times K_{1j}^I(s)\right)
 \end{array}\right\}
 \mbox{  is robustly Hurwitz stable.  }
\end{array}
$$

{\bf Theorem 2}\quad ${\cal G}(s)$ is robustly Hurwitz stable if and only if
${\cal G}_E(s)$ is robustly Hurwitz stable.

{\bf Remark:}\ \ Theorem 2 is consistent with the result in \cite{Mori}. In \cite{Mori}, the authors
obtained their result using some theorem in signal processing. Our proof is based on the properties of
matrix determinant, hence is more straightforward and self-contained.

\section{Example}

Consider the uncertain polynomial matrix
$$
\begin{array}{l}
{\cal A}(s)=\left\{(a_{ij}(s))_{3\times 3}\right\}\\
a_{11}(s)=\lambda_{11}b1_{11}(s)+(1-\lambda_{11})b2_{11}(s)\\
a_{12}(s)=\lambda_{12}b1_{12}(s)+(1-\lambda_{12})b2_{12}(s)\\
a_{13}(s)=\lambda_{13}b1_{13}(s)+(1-\lambda_{13})b2_{13}(s)\\
a_{21}(s)=\lambda_{21}b1_{21}(s)+(1-\lambda_{21})b2_{21}(s)\\
a_{22}(s)=\lambda_{22}b1_{22}(s)+(1-\lambda_{22})b2_{22}(s)\\
a_{23}(s)=\lambda_{23}b1_{23}(s)+(1-\lambda_{23})b2_{23}(s)\\
a_{31}(s)=\lambda_{31}b1_{31}(s)+(1-\lambda_{31})b2_{31}(s)\\
a_{32}(s)=\lambda_{32}b1_{32}(s)+(1-\lambda_{32})b2_{32}(s)\\
a_{33}(s)=\lambda_{33}b1_{33}(s)+(1-\lambda_{33})b2_{33}(s)
\end{array}
$$
Then,
$$
\begin{array}{l}
E_{ij}(s)=\{ a_{ij}(s)\}\\
K_{ij}(s)=\{b1_{ij}(s),b2_{ij}(s) \}\\
S_3=\{\sigma_1,\dots,\sigma_6\}\\
\sigma_1: 1\rightarrow 1; 2\rightarrow 2;3\rightarrow 3\\
\sigma_2: 1\rightarrow 2; 2\rightarrow 3;3\rightarrow 1\\
\sigma_3: 1\rightarrow 3; 2\rightarrow 1;3\rightarrow 2\\
\sigma_4: 1\rightarrow 1; 2\rightarrow 3;3\rightarrow 2\\
\sigma_5: 1\rightarrow 2; 2\rightarrow 1;3\rightarrow 3\\
\sigma_6: 1\rightarrow 3; 2\rightarrow 2;3\rightarrow 1
\end{array}
$$
Let
$$
{\cal A}_E(s)=\bigcup_{\sigma\in S_3 }\left \{(a_{ij}(s))_{3 \times 3}:
a_{ij}(s)\left\{
\begin{array}{l}
\in E_{ij}(s)  \mbox{  if   } i=\sigma(j) \\
 \in K_{ij}(s) \mbox{  if   } i\not=\sigma(j)
\end{array}\right.
\right \}
$$
By Theorem 1, ${\cal A}(s)$ is robustly $D$ stable if and only if
${\cal A}_E(s)$ is robustly $D$ stable.

\section{Conclusions}

This paper discussed the robust $D$-stability problems for MIMO uncertain systems.
The Edge Theorem and Kharitonov Theorem have been generalized to multivariable case.

\end{document}